\documentclass{article}

\usepackage{amsmath,amssymb,amsthm}

\newtheorem{theorem}{Theorem}[section]
\newtheorem{corollary}[theorem]{Corollary}
\newtheorem{proposition}[theorem]{Proposition}

\theoremstyle{remark}
\newtheorem{remark}{Remark}[section]
\theoremstyle{definition}
\newtheorem{definition}{Definition}[section]

\theoremstyle{definition}
\newtheorem{example}{Example}[section]
\newtheorem{conjecture}[theorem]{Conjecture}
\begin{document}

\title{Carath\'{e}odory-Equivalence, Noether Theorems, 
and Tonelli Full-Regularity in the Calculus of Variations
and Optimal Control\footnote{The contents of this paper are taken
from the author's Ph.D. thesis, University of Aveiro, 2002.
(Supervisor: A.~V. Sarychev). The paper is
submitted for publication in a Special Issue of the 
\emph{J. of Mathematical Sciences}. To be presented
at the \emph{First Junior European Meeting
``Control Theory and Stabilization''}, Dijon, October 2-4 2002.
The date of this version is June 21, 2002.}}

\author{Delfim F. M. Torres\\
        \texttt{delfim@mat.ua.pt}}

\date{Department of Mathematics\\
      University of Aveiro\\
      3810-193 Aveiro, Portugal\\
      \texttt{http://www.mat.ua.pt/delfim}}

\maketitle

%%%%%%%%%%%%%%%%%%%%%%%%%%%

\begin{abstract}
We study, in a unified way, the following
questions related to the properties of Pontryagin extremals for 
optimal control problems with unrestricted controls:
i) How the transformations, which define the equivalence of 
two problems, transform the extremals?
ii) How to obtain quantities which are conserved along any extremal?
iii) How to assure that the set of extremals include the minimizers predicted 
by the existence theory? 
These questions are connected to: 
i) the Carath\'{e}odory method which establishes a correspondence 
between the minimizing curves of equivalent problems;
ii) the interplay between the
concept of invariance and the theory of optimality conditions
in optimal control, which are the concern of the theorems of Noether;
iii) regularity conditions for the minimizers
and the work pioneered by Tonelli.
\end{abstract}

%%%%%%%%%%%%%%%%%%%%%%%%%%%

\section{Introduction}

For more than three centuries, the calculus of variations 
played a central role stimulating the development of mathematics
and the development of physics. Today, the calculus
of variations, and its natural generalization known as the theory
of optimal control, remain relevant and useful, generating new exciting 
and deep questions. There is a substantial progress in fundamental issues
of both theory and applications \cite{MR1796845,300Sussmann}.
In this paper we will address some of these questions.

We study properties of the extremals and minimizers for various problems
of the calculus of variations and optimal control.
We are particularly interested in problems which may appear quite
different but still can be reduced to the same problem if one
uses appropriate transformations. Such problems are said to be equivalent.
It turns out that for the equivalent problems there is a direct relation between 
admissible state-control pairs and the value for the cost functionals. 
In particular, after \emph{solving} one problem, it is then straightforward
to obtain the solutions for all the equivalent problems
from the transformations which define the equivalence.

The standard scheme to solve
a problem in the calculus of variations or optimal control 
proceeds along the following three steps.
First we prove that a solution to the problem exists. 
Second we assure the applicability of necessary optimality conditions.
Finally we apply the necessary conditions which identify
the extremals (the candidates). Further elimination, if necessary,
identifies the minimizer(s) of the problem.
As pointed out by L.~C.~Young \cite{MR41:4337},
although both the calculus of variations and optimal
control have born from the study of necessary optimality
conditions, any such theory is ``naive'' 
until the existence of minimizers is assured.
The process leading to existence theorems was introduced
by Leonida Tonelli, in the years 1911-1915, through the
so called \emph{direct method}.
It turns out that, even for the simplest problem of the calculus
of variations, the hypotheses of the existence
theory do not imply those of the necessary optimality
conditions. This is to say that
all the three steps in the above procedure are indeed crucial:
it does not make sense to apply necessary optimality conditions
if no solution to the problem exists;
and it may be the case that a solution exist but fails to satisfy the standard 
necessary optimality conditions such as the Euler-Lagrange equations 
or the Pontryagin maximum principle. Therefore, regularity conditions
are also an essential step in the process of solving a problem in the calculus of
variations or optimal control. They close the gap between existence and
optimality theories, assuring that all the minimizers are indeed extremals
\cite{MR2001j:49062}.

The study of equivalent problems in the calculus of variations 
and optimal control is not enough. It is also important to know
how the extremals of the problems are related. 
In the terminology of Constantin Carath\'{e}odory \cite[\S 227]{ZBL0505.49001},
two problems of the calculus of variations are said to be equivalent
if the respective Lagrangians differ by a total derivative. The
importance of this equivalence concept is due to the fact
that it implies that the Euler-Lagrange equations are identical 
for both problems. One can say that for Carath\'{e}odory
is the correspondence between the extremals,
and not that of the problems, the key concept to define equivalence.
To the best of our knowledge this concept of Carath\'{e}odory-equivalence 
has not been previously explored, or even considered, 
in the optimal control context.
Here we will be mainly interested in the following trivial 
but important remark: two Carath\'{e}odory-equivalent problems
have ``the same'' \emph{conservation laws}. This is not necessarily the case
for equivalent problems: equivalence does not imply 
Carath\'{e}odory-equivalence and the other way around. 
Surprising enough, when one restrict attention to the abnormal extremals,
the two concepts seem to be quite the same.
\begin{conjecture}
\label{Torres:My:Conjecture}
Two problems of optimal control are equivalent if, and only if,
they are abnormal-Carath\'{e}odory-equivalent.
\end{conjecture}
We will show in Section~\ref{Torres:Sec:CE} the validity of 
Conjecture~\ref{Torres:My:Conjecture} for
equivalent problems under transformations of the type of Gamkrelidze
\cite[\S 8.5]{MR58:33350c} and under a time-reparameterization
introduced by the author in \cite{delfimEJC}.

Conservation laws, that is, conserved quantities along the 
extremals of the problem, are obtained in the calculus of variations 
with the help of the famous symmetry theorems of Emmy Noether.
These classical results are known as the (first) Noether 
theorem and the second Noether theorem, and 
explain the correspondence between invariance of the problem with respect
to a family of transformations and the existence of conservation laws.
The first Noether theorem establishes the existence of $\rho$
conservation laws of the Euler-Lagrange differential equations when
the Lagrangian $L$ is invariant under a family of
transformations containing $\rho$ parameters.
The second Noether theorem establishes the existence of $k \left(m + 1\right)$
conservation laws when the Lagrangian is invariant under 
a family of transformations which, rather than dependence on
parameters, as in the first theorem,
depend upon $k$ arbitrary functions and their derivatives up to
order $m$. Extensions of the first theorem for the Pontryagin extremals of optimal
control problems are available in \cite{delfimEJC,delfim3ncnw,torresMED2002}.
In Section~\ref{Torres:Sec:NT} we provide a rather general formulation
of the first Noether theorem which envolves all the peculiarities of previous
results. For optimal control versions of the second theorem we refer the reader 
to \cite{torresCM02I06E}.
We will argue in Section~\ref{Torres:Sec:TFR} 
that the conservation laws obtained from
the use of Noether's first theorem play an important role in the
acquisition of regularity conditions. 

Further extensions and related results are possible.
Due to the restrictions on the volume of the paper
we can not provide them here. 
The reader can find more details, complete and detailed proofs,
illustrative examples, additional material and a complete list of references,
in the author's thesis \cite{torresPhD}, available in Portuguese.

%%%%%%%%%%%%%%%%%%%%%%%%%%%

\section{Carath\'{e}odory-Equivalence}
\label{Torres:Sec:CE}

In this section we address, under two types of transformations
involving change of the time-variable $t$, 
the following question:
\emph{How the transformations, which define the equivalence of 
two problems, affect the Pontryagin extremals?} First we
will need to introduce the problem considered in optimal 
control theory, and to give a characterization
of the Pontryagin extremals.

The objective of the paper is to study some properties of the minimizing
trajectories for general problems of optimal control
in the case where controls are unconstrained, like
in the calculus of variations. We will be considering,
without any loss of generality, problems in the Lagrange
form. This is indeed a general problem, and 
Bolza type problems or Mayer type problems can be put
easily in this form. 
We look for a pair $\left(x(\cdot),u(\cdot)\right)$,
satisfying a control dynamical equation
described by a system of ordinary differential equations
\begin{equation}
\label{Torres:eq:controlEquation}
\dot{x}(t) = \varphi\left(t,x(t),u(t)\right) \, ,
\end{equation}
in such a way the pair $\left(x(\cdot),u(\cdot)\right)$ minimizes a given 
integral functional:
\begin{equation*}
I\left[x(\cdot),u(\cdot)\right] 
= \int_{a}^{b}L\left(t,x(t),u(t)\right) \mathrm{d}t \longrightarrow \min \, .
\end{equation*}
This problem is denoted by $(P)$.
The state trajectories $x(\cdot)$ are assumed to be absolute
continuous functions and the admissible controls $u(\cdot)$
to be Lebesgue integrable: 
\begin{equation*}
x(\cdot) \in W_{1,1}\left(\left[a,b\right];\,\mathbb{R}^{n}\right)\, , \quad
u(\cdot) \in L_{1}\left(\left[a,b\right];\,\mathbb{R}^{r}\right) \, .
\end{equation*}
For simplicity of exposition, we will assume that the
Lagrangian $L$ and function $\varphi$ are $C^1$-smooth with respect
to all variables. The results of the paper
are valid for all kinds of boundary conditions one
may want to consider. For this reason, boundary conditions are not
considered in our formulation of the optimal control problem.
We remark that, \emph{a priori}, optimal controls may be unbounded and that
the problems of the calculus of variations, like the basic problem
of the calculus of variations or the problems with high-order derivatives,
can be reduced in the obvious way to the Lagrange problem $(P)$.

Both the calculus of variations and optimal
control theory have born from the study of first-order necessary optimality
conditions: the Euler-Lagrange equations,
in the case of the calculus of variations, 
which appear in the Euler's celebrated monograph of
1744, and the Pontryagin maximum principle, 
in the case of optimal control,
which appear $\sim 1956$.
The Pontryagin maximum principle gives conditions
under which all the minimizers are Pontryagin extremals.
\begin{definition}
The quadruple $\left(x(\cdot),u(\cdot),\psi_{0},\psi(\cdot)\right)$, 
$\psi_0 \le 0$, $\psi\left(\cdot\right) \in W_{1,1}$,
$(\psi_0,\,\psi(\cdot)) \ne 0$,
is a Pontryagin extremal if it satisfies:
\begin{description}
\item[] the Hamiltonian system
\begin{equation*}
\dot{x}=\dfrac{\partial H}{\partial\psi} \, , \quad
\dot{\psi}= - \dfrac{\partial H}{\partial x} \, ;
\end{equation*}
\item[] the maximality condition
\begin{equation*}
H\left(t,x(t),u(t),\psi_{0},\psi(t)\right)
=\max_{u\in\mathbb{R}^{r}} H\left(t,x(t),u,\psi_{0},\psi(t)\right) \, ;
\end{equation*}
\end{description}
with the Hamiltonian 
$H(t,x,u,\psi_0,\psi)=\psi_{0}\,L\left(t,x,u\right) + \psi \cdot \varphi\left(t,x,u\right)$.
An extremal is said to be abnormal when $\psi_0$ vanishes
and normal otherwise. 
\end{definition}
The first equation in the Hamiltonian system is just the control equation
\eqref{Torres:eq:controlEquation}. The second is known as the \emph{adjoint system}.

For the basic problem of the calculus of variations one
has $\dot{x} = u$ and the Hamiltonian is given by $H = \psi_0 L + \psi \cdot u$.
From the adjoint system we obtain 
\begin{equation}
\label{Torres:AdjSys:PBCV}
\dot{\psi} = - \psi_0 \dfrac{\partial L}{\partial x} \, ,
\end{equation}
while from the maximality condition one gets 
\begin{equation}
\label{Torres:MaxCond:PBCV}
\psi = - \psi_0 \dfrac{\partial L}{\partial u} \, .
\end{equation}
The maximum principle asserts that $\psi_0$ and $\psi(\cdot)$ do not vanish 
simultaneously and it comes immediately, from \eqref{Torres:MaxCond:PBCV},
that no abnormal minimizers exist for the basic problem of the calculus of variations.
From equalities \eqref{Torres:AdjSys:PBCV} and \eqref{Torres:MaxCond:PBCV}
one concludes that if $x(\cdot)$ is a minimizer
it satisfies the Euler-Lagrange equations:
\begin{equation}
\label{Torres:eq:EL}
\frac{\mathrm{d}}{\mathrm{d}t} \frac{\partial L}{\partial u}\left(t,x(t),\dot{x}(t)\right) 
= \frac{\partial L}{\partial x}\left(t,x(t),\dot{x}(t)\right) \, .
\end{equation}
A solution $x(\cdot)$ of \eqref{Torres:eq:EL} is called an 
(Euler-Lagrange) extremal.
The classical conditions \cite{MR29:3316b}, assuring that
all minimizers are extremals, assume that the Lagrangian
$L$ and function $\varphi$ are continuous with respect
to all variables and continuously differentiable with
respect to the state variables $x$; while the optimal
controls are assumed to be essentially bounded: 
$L(\cdot,x,\cdot),\,\varphi(\cdot,x,\cdot) \in C$,
$L(t,\cdot,u),\,\varphi(t,\cdot,u) \in C^1$,
$u(\cdot) \in L_{\infty}$.
For the basic problem of the calculus of variations, this means
that the Euler-Lagrange equations are valid for minimizers in the
class of Lipschitzian functions.
Proving general versions of the maximum
principle under weaker hypotheses is still very much
in progress \cite{MR1806192}. 
Conditions with $u(\cdot) \in L_{1}$ do exist, but
they postulate growth conditions on the Lagrangian $L$
and functions $\varphi$. For example, the following conditions follow easily  
from Berkovitz's \cite{MR51:8914} or 
Clarke's version of the maximum principle \cite{MR54:3540}:
\begin{equation}
\label{Torres:CondApplPMPBerkClar}
\left\|\frac{\partial L}{\partial x}\right\| \le c \left|L\right| + k \, , \quad
\left\|\frac{\partial \varphi_i}{\partial x}\right\| \le c \left|\varphi_i\right| + k \, ,
\end{equation}
$i = 1,\ldots,n$, for some constants $c>0$ and $k$.

The following theorem provides an interesting 
property of the Pontryagin extremals.
\begin{theorem}[\cite{delfimIO}]
\label{Torres:Cap3:r:mainresult}
Let $F\left(t,x,u,\psi_0,\psi\right)$,
$F : [a,b] \times \mathbb{R}^n \times \mathbb{R}^r \times \mathbb{R}_{0}^{-}
\times \mathbb{R}^n \rightarrow \mathbb{R}$, be a continuous differentiable 
function with respect to $t$, $x$, $\psi_0$, and $\psi$, for fixed $u$. If there exists
a function $G(\cdot) \in L_1\left([a,b] ; \mathbb{R} \right)$ 
such that
\begin{equation*}
\left\|\nabla_{(t,x,\psi)}
F\left(t,x(t),u(s),\psi_0,\psi(t)\right)\right\| \le G(t) \, 
\quad \left(s, t \in [a,b]\right) \, ,
\end{equation*}
and for almost all $t$'s in the interval $[a,b]$ the condition
\begin{equation*}
F\left(t,x(t),u(t),\psi_0,\psi(t)\right)
= \max_{v \in \mathbb{R}^r} F\left(t,x(t),v,\psi_0,\psi(t)\right)
\end{equation*}
is true along the Pontryagin extremals
$\left(x(\cdot),u(\cdot),\psi_0,\psi(\cdot)\right)$ of the optimal
control problem $(P)$, then 
$t \rightarrow F\left(t,x(t),u(t),\psi_0,\psi(t)\right)$ is
absolutely continuous and the equality
\begin{equation*}
\frac{\mathrm{d}F}{\mathrm{d}t} = \frac{\partial F}{\partial t}
+ \frac{\partial F}{\partial x} \cdot \frac{\partial H}{\partial \psi}
- \frac{\partial F}{\partial \psi} \cdot \frac{\partial H}{\partial x}
\doteq \frac{\partial F}{\partial t} + \left\{F,H\right\} 
\end{equation*}
holds along the extremals,
where $\left\{F,H\right\}$ denotes the Poisson bracket
of the functions $F$ and $H$,
and on the left-hand side we have the total derivative with respect
to $t$, and on the right-hand side partial derivatives.
\end{theorem}
If one chooses $F$ in Theorem~\ref{Torres:Cap3:r:mainresult}
to be the Hamiltonian $H$, one gets: 
\begin{theorem}
If $\left(x(\cdot),u(\cdot),\psi_0,\psi(\cdot)\right)$ is
a Pontryagin extremal, then the function 
$H\left(t,x(t),u(t),\psi_0,\psi(t)\right)$
is an absolutely continuous function of $t$ and satisfies the equality
\begin{equation}
\label{Torres:eq:dHdtpartialHdt}
\frac{\mathrm{d}H}{\mathrm{d}t}\left(t,x(t),u(t),\psi_0,\psi(t)\right)
= \frac{\partial H}{\partial t}\left(t,x(t),u(t),\psi_0,\psi(t)\right) \, .
\end{equation}
\end{theorem}
Equation \eqref{Torres:eq:dHdtpartialHdt}
corresponds, for the basic problem of the calculus of variations,
to the classical DuBois-Reymond necessary condition:
\begin{equation}
\label{Torres:eq:DuBoisReymond}
\frac{\mathrm{d}}{\mathrm{d}t} \left[
L\left(t,x(t),\dot{x}(t)\right) - \dot{x}(t) \cdot
\frac{\partial L}{\partial u}\left(t,x(t),\dot{x}(t)\right)\right]
= \frac{\partial L}{\partial t}\left(t,x(t),\dot{x}(t)\right) \, .
\end{equation}

From our Theorem~\ref{Torres:Cap3:r:mainresult},
a necessary and sufficient condition for a function
$F$ to be a conservation law is immediately obtained.
\begin{corollary}
\label{Torres:CM:RemProp:APDIO}
Under the conditions of Theorem~\ref{Torres:Cap3:r:mainresult},
$F(t,x,u,\psi_0,\psi)$ is constant along every
Pontryagin extremal of the problem if, and only if,
\begin{equation*}
\frac{\partial F}{\partial t}
+ \frac{\partial F}{\partial x} \cdot \frac{\partial H}{\partial \psi}
- \frac{\partial F}{\partial \psi} \cdot \frac{\partial H}{\partial x}
= 0 \, .
\end{equation*}
\end{corollary}
Corollary~\ref{Torres:CM:RemProp:APDIO} is very
useful for the characterization of optimal control problems
with given conserved quantities along the Pontryagin extremals.
For example, if one wants to find a problem for which the function 
\begin{equation}
\label{Torres:eq:CLCh1}
F = H \psi x
\end{equation}
is constant in $t$ along the respective
extremals, a necessary and sufficient condition is given by the relation
\begin{equation*}
\psi x \frac{\partial H}{\partial t}
+ \psi H \frac{\partial H}{\partial \psi}
- H x \frac{\partial H}{\partial x} = 0 \, .
\end{equation*}
One such problem is therefore
\begin{gather}
\label{Torres:LuFiuvezx}
\int_a^b L\left(u(t)\right) \mathrm{d}t \longrightarrow \min \, , \\  
\dot{x}(t) = \varphi\left(u(t)\right) x(t) \notag \, .
\end{gather}
Usually, the form of the problem is already known and
conditions are sought in such a way that the problem as
some ``good'' properties. Let us see one such situation.
The following problem is related to the study of 
cubic polynomials on Riemannian manifolds:
\begin{gather}
\int_0^T \sum_{i=1}^{n} \left(u_i(t)\right)^2 \mathrm{d}t 
\longrightarrow \min \, , \notag \\
\label{Torres:eq:pr:cubicPolyRiemMani}
\begin{cases}
\dot{x_1}(t) = x_2(t) \, , \\
\dot{x_2}(t) = \sum_{i=1}^{n} X_i\left(x_1(t)\right) u_i(t)\, .
\end{cases} 
\end{gather}
The problem is autonomous and from \eqref{Torres:eq:dHdtpartialHdt}
we know that the respective Hamiltonian $H$ is conserved along the extremals.
Determination of the explicit solutions to problem 
\eqref{Torres:eq:pr:cubicPolyRiemMani}
is a difficult task and is, in general, an open problem.
However, the extremals can be explicitly computed if a new independent
conserved quantity is found. The question is:
what kind of conditions shall we impose
on the vector fields $X_i$ in order to obtain the new conserved quantity?
To answer this question one needs to solve a characterization problem as before.
Let $F = k_1 \psi_1 x_1 + k_2 \psi_2 x_2$
($k_1$ and $k_2$ constants). This is called in the literature
a \emph{momentum map}. Using the relation given by 
Corollary~\ref{Torres:CM:RemProp:APDIO} one gets
\begin{multline*}
k_1 \psi_1 x_2 + k_2 \psi_2 \left(X_1(x_1) u_1 + \cdots + X_n(x_1) u_n\right) \\
- k_1 x_1 \psi_2 \left(X'_1(x_1) u_1 + \cdots + X'_n(x_1) u_n\right)
- k_2 x_2 \psi_1 = 0 \, .
\end{multline*}
This condition is trivially satisfied 
if $k_1 = k_2$ and $X'_i(x_1) x_1 = X_i(x_1)$, $i = 1,\,\ldots,\,n$.
We have just proved the following proposition.
\begin{proposition}
\label{Torres:prop:homogenei}
If the homogeneity condition 
\begin{equation}
\label{Torres:eq:homogeneidade}
X_i\left(\lambda x_1\right) = \lambda X_i(x_1) \, , \quad
i = 1,\,\ldots,\,n \, , \, \forall \, \lambda > 0 \, , 
\end{equation}
holds, then
\begin{equation*}
\psi_1(t) x_1(t) + \psi_2(t) x_2(t)
\end{equation*}
is constant in $t \in [0,T]$ along any extremal of the problem
\eqref{Torres:eq:pr:cubicPolyRiemMani}.
\end{proposition}
We shall elaborate more on this issue later, in Section~\ref{Torres:Sec:NT},
in relation with Noether's theorems (\textrm{cf.} Example~\ref{Torres:ex:34}). 

After this short introduction to Pontryagin extremals and
their characterization, we are now in conditions to study
how the extremals are affected when one transforms problem $(P)$.
Let us consider the following optimal control problem:
\begin{gather}
J\left[t(\cdot),z(\cdot),v(\cdot)\right] = 
\int_{\tau_a}^{T} \Upsilon\left(t(\tau),z(\tau),v(\tau)\right)
L\left(t(\tau),z(\tau),v(\tau)\right)
\mathrm{d}\tau \, \longrightarrow \min \, , \notag \\
\begin{cases}
t'(\tau) &= \Upsilon\left(t(\tau),z(\tau),v(\tau)\right) \\
z'(\tau) &= \Upsilon\left(t(\tau),z(\tau),v(\tau)\right) 
\varphi\left(t(\tau),z(\tau),v(\tau)\right)
\end{cases} \label{Torres:eq:MudVarIndCG} \\
v : \mathbb{R} \rightarrow \mathbb{R}^r \, , \notag \\
t(\tau_a) = a \, , \quad t(T) = b \, , \notag 
\end{gather}
where $\Upsilon(\cdot,\cdot,\cdot)$ is a strictly positive 
continuously differentiable function,
\begin{equation*}
C^1 \ni \Upsilon(t,z,v) : \mathbb{R} \times \mathbb{R}^{n} \times 
\mathbb{R}^{r} \rightarrow \mathbb{R}^{+} \, ,
\end{equation*}
and $T$ is free. Problem \eqref{Torres:eq:MudVarIndCG} is obtained
from $(P)$ by entering a new time variable $\tau$, related with with
$t$ by the relation
\begin{equation}
\label{Torres:eq:MyOlk44}
\tau(t) = \tau_a + \int_a^t \frac{1}{\Upsilon\left(\theta,x(\theta),u(\theta)\right)} 
\mathrm{d} \theta \, ,
\quad t \in [a,b] \, .
\end{equation}
Compared with $(P)$, problem \eqref{Torres:eq:MudVarIndCG} 
has one more state variable. Namely, 
its state variables are $t(\cdot)$ and $z(\cdot)$.
We note that the problem is autonomous:
both the Lagrangian and the right-hand side of the control system
do not depend directly on $\tau$. Thereafter, the admissible set of problem \eqref{Torres:eq:MudVarIndCG} is invariant with respect to translations 
on the time variable $\tau$. For the concrete situation wherein $\tau_a = 0$ and
$\Upsilon(t,z,v) = \frac{1}{L\left(t,z,v\right)}$,
one obtains the transformation introduced by 
R.~V.~Gamkrelidze \cite[Chap. 8]{MR58:33350c},
of the Lagrange problem $(P)$ into the autonomous time optimal problem. 

We denote by $H$ the Hamiltonian associated with problem $(P)$,
and by $\mathcal{H}$ the Hamiltonian associated with \eqref{Torres:eq:MudVarIndCG}:
\begin{gather*}
H\left(t,x,u,\psi_0,\psi\right) = \psi_0 L(t,x,u) + \psi \cdot \varphi(t,x,u) \, , \\
\mathcal{H}\left(t,z,v,p_0,p_t,p_z\right) = 
\left[p_0 L(t,z,v) + p_t + p_z \cdot \varphi(t,z,v)\right] \Upsilon(t,z,v) \, .
\end{gather*}
As far as $\mathcal{H}$ does not depend on $\tau$, it follows
from Theorem~\ref{Torres:Cap3:r:mainresult} that
$\mathcal{H} \equiv const$ along an extremal. 
The following theorems assert that 
the extremals of problem $(P)$ are related to those extremals of problem
\eqref{Torres:eq:MudVarIndCG} for which this constant is zero.
\begin{theorem}[\cite{torresPhD}]
\label{Torres:Cap4:T1NP}
Let $\left(t(\tau),z(\tau),v(\tau),p_0,p_t(\tau),p_z(\tau)\right)$,
$\tau \in [\tau_a,T]$, be a Pontryagin extremal of problem
\eqref{Torres:eq:MudVarIndCG} with
\begin{equation}
\label{Torres:Cap4:ZeroLevel}
\mathcal{H}\left(t(\tau),z(\tau),v(\tau),p_0,p_t(\tau),p_z(\tau)\right) = 0 \, .
\end{equation}
Then
$\left(x(t),u(t),\psi_0,\psi(t)\right) = 
\left(z(\tau(t)),v(\tau(t)),p_0,p_z(\tau(t))\right)$,
$t \in [a,b]$, 
where $\tau(\cdot)$ is the inverse function of $t(\cdot)$,
is a Pontryagin extremal of $(P)$. Moreover, the value
for the functionals coincide:
\begin{equation}
\label{Torres:P:PGam:MVF}
I\left[x(\cdot),u(\cdot)\right] = J\left[t(\cdot),z(\cdot),v(\cdot)\right] \, .
\end{equation}
\end{theorem}

\begin{theorem}[\cite{torresPhD}]
\label{Torres:Cap4:ZeroLevel:OWA}
Let $\left(x(t),u(t),\psi_0,\psi(t)\right)$,
$t \in [a,b]$, be a Pontryagin extremal of problem $(P)$.
Then, with $t(\cdot)$ the inverse function of \eqref{Torres:eq:MyOlk44},
\begin{multline*}
\left(t(\tau),z(\tau),v(\tau),p_0,p_t(\tau),p_z(\tau)\right)
= \left(t(\tau),x(t(\tau)),u(t(\tau)),\psi_0, \right. \\
\left. -H\left(t(\tau),x(t(\tau)),u(t(\tau)),\psi_0,%
\psi(t(\tau))\right),\psi(t(\tau))\right) \, ,
\end{multline*}
$\tau \in [\tau_a,T]$, is a Pontryagin extremal of \eqref{Torres:eq:MudVarIndCG}
which satisfies equalities \eqref{Torres:Cap4:ZeroLevel} and \eqref{Torres:P:PGam:MVF}.
\end{theorem}
We remark that the correspondence given by Theorems~\ref{Torres:Cap4:T1NP}
and \ref{Torres:Cap4:ZeroLevel:OWA} keeps the normality or abnormality
of the extremals ($\psi_0 = p_0$). 

Similar correspondences between the extremals, to the one
established by Theorems~\ref{Torres:Cap4:T1NP}
and \ref{Torres:Cap4:ZeroLevel:OWA}, can be obtained
under different transformations of the problem $(P)$.
One such transformation
is based on an idea of time reparameterization
introduced by the author in \cite{delfimEJC}.
The idea generalizes a well known time reparameterization that 
has proved to be useful in many different contexts of the calculus 
of variations and optimal control (see references in \cite{delfimEJC}).
Considering $t$ as a dependent variable, we introduce
a one to one  Lipschitzian transformation 
$[a,b] \ni t \longmapsto \tau \in [a,b]$,
$\frac{\mathrm{d}}{\mathrm{d} \tau} t(\tau) > 0$, 
such that
\begin{gather*}
L\left(t,x(t),u(t)\right)\,\mathrm{d}t =
L\left(t(\tau),x(t(\tau)),u(t(\tau))\right)\,
\frac{\mathrm{d}t(\tau)}{\mathrm{d}\tau}\,\mathrm{d} \tau \, , \\
\frac{\mathrm{d}}{\mathrm{d} \tau} x\left(t(\tau)\right) =
\frac{\mathrm{d} x\left(t(\tau)\right)}{\mathrm{d} t}\,\frac{\mathrm{d} t(\tau)}{\mathrm{d} \tau}
= \varphi\left(t(\tau),x\left(t(\tau)\right),u\left(t(\tau)\right)\right)\,
\frac{\mathrm{d} t(\tau)}{\mathrm{d}\tau} \, .
\end{gather*}
In this way, if one consider the notations
$z(\tau) = x(t(\tau))$ and $w(\tau) = u(t(\tau))$, problem $(P)$ 
takes the form
\begin{equation}
\label{Torres:eq:Ptau}
\begin{gathered}
K\left[t(\cdot),z(\cdot),v(\cdot),w(\cdot)\right] = 
\int_{a}^{b}
L\left(t(\tau),z(\tau),w(\tau)\right)\,v(\tau)\,\mathrm{d} \tau
\longrightarrow \min \\ \tag{$P_{\tau}$}
\begin{cases}
t'(\tau) = v(\tau) \\
z'(\tau) = \varphi\left(t(\tau),z(\tau),w(\tau)\right)\,v(\tau)
\end{cases}
\end{gathered}
\end{equation}
\begin{equation*}
\begin{gathered}
t(a) = a \, , \quad t(b) = b \, , \\
t(\cdot) \in W_{1,\infty}\left([a,b];\,[a,b]\right)\, , \quad
z(\cdot) \in W_{1,1}\left([a,b];\,\mathbb{R}^n\right)\, , \\
v(\cdot) \in L_{\infty}\left([a,b];\,\mathbb{R}^{+}\right)\, , \quad
w(\cdot) \in L_{1}\left([a,b];\,\mathbb{R}^r\right) \, .
\end{gathered}
\end{equation*}
For the new transformed problem \eqref{Torres:eq:Ptau},
the state variables are $t(\tau)$ and $z(\tau)$ while 
the controls are $v(\tau)$ and $w(\tau)$. The fact that the control
variable $v(\cdot)$ takes only strictly positive values,  
assure that $t(\tau)$ has an inverse function $\tau(t)$.

Next theorem shows how to construct an extremal of $(P)$
given an extremal of problem \eqref{Torres:eq:Ptau}.

\begin{theorem}[\cite{torresMCSS}]
\label{Torres:r:propImpAbn}
Let $\left(t(\cdot),z(\cdot),v(\cdot),w(\cdot),p_0,p_t(\cdot),p_z(\cdot)\right)$ 
be an extremal of \eqref{Torres:eq:Ptau}. Then
\begin{equation*}
\left(x(\cdot),u(\cdot),\psi_0,\psi(\cdot)\right) =
\left(z\left(\tau(\cdot)\right),w\left(\tau(\cdot)\right),
p_0,p_z\left(\tau(\cdot)\right)\right) \, ,
\end{equation*}
where $\tau(\cdot)$ is the inverse function of $t(\cdot)$,
is an extremal of $(P)$ with the same value for
the functional: $I\left[x(\cdot),u(\cdot)\right] = 
K\left[t(\cdot),z(\cdot),v(\cdot),w(\cdot)\right]$.
\end{theorem}

The transformed problem \eqref{Torres:eq:Ptau} is autonomous, and we already
know that the corresponding Hamiltonian is constant along the extremals.
As before, for the Gamkrelidze-type transformations,
to a Pontryagin extremal of the original problem $(P)$ one can correspond
an extremal of the transformed problem for which the Hamiltonian
vanishes.
\begin{theorem}[\cite{torresMCSS}]
\label{Torres:Cap4:r:prop2}
Let $\left(x(\cdot),u(\cdot),\psi_0,\psi(\cdot)\right)$
be a Pontryagin extremal of $(P)$. Then, for all
$v(\cdot) \in L_{\infty}\left([a,b];\,\mathbb{R}^{+}\right)$ such that
$\int_{a}^{b} v(\theta)\,\mathrm{d}\theta = b - a$, the 7-uple
$\left(t(\cdot),z(\cdot),v(\cdot),w(\cdot),p_0,p_t(\cdot),p_z(\cdot)\right)$, 
defined by
\begin{gather*}
t(\tau) = a + \int_{a}^{\tau} v(\theta)\,\mathrm{d}\theta\, ,  \\
z(\tau) = x\left(t(\tau)\right)\, , \quad w(\tau) = u\left(t(\tau)\right)\, , \\
p_0 = \psi_0\, , \quad p_z(\tau) = \psi(t(\tau))\, ,  \\
p_t(\tau) =-H\left(t(\tau),x(t(\tau)),u(t(\tau)),\psi_0,\psi(t(\tau))\right) \, ,
\end{gather*}
is a Pontryagin extremal of \eqref{Torres:eq:Ptau}
giving a zero value for the respective Hamiltonian. Moreover,
$K\left[t(\cdot),z(\cdot),v(\cdot),w(\cdot)\right] =
I\left[x(\cdot),u(\cdot)\right]$.
\end{theorem}

From Theorems~\ref{Torres:r:propImpAbn} and \ref{Torres:Cap4:r:prop2}
the following corollary is trivially obtained. 
\begin{corollary}
If $\left(\tilde{x}(t),\tilde{u}(t)\right)$ 
is a minimizer of $(P)$, then the quadruple
\begin{equation*}
\left(\tilde{t}(\tau),\tilde{z}(\tau),\tilde{v}(\tau),\tilde{w}(\tau)\right)
=(\tau,\tilde{x}(\tau),1,\tilde{u}(\tau)) 
\end{equation*}
furnishes a minimizer to \eqref{Torres:eq:Ptau}.
\end{corollary}

Once again, Theorems~\ref{Torres:r:propImpAbn} and \ref{Torres:Cap4:r:prop2}
establish a correspondence between the abnormal extremals of the original
and transformed problem.
\begin{corollary}
If no abnormal extremals exist for the problem $(P)$, 
then no abnormal extremals exist for the problem \eqref{Torres:eq:Ptau} too. 
If no abnormal extremals exist for the problem \eqref{Torres:eq:Ptau}, 
then no abnormal extremals exist for the problem $(P)$ too. 
\end{corollary}

We shall see in Section~\ref{Torres:Sec:TFR} that one
can obtain regularity conditions,
assuring that all minimizing controls of $(P)$, predicted by
Tonelli's existence theorem, are Pontryagin extremals, 
from the applicability conditions
of the Pontryagin maximum principle to
the transformed problems. The proof relies on certain
conserved quantities along the Pontryagin extremals.
These conserved quantities are addressed in the following section.

%%%%%%%%%%%%%%%%%%%%%%%%%%%

\section{Noether Symmetry Theorems}
\label{Torres:Sec:NT}

We now turn our attention to the following question:
\emph{How to obtain quantities which are conserved along the 
Pontryagin extremals?}
This is an important, profound, and far-reaching litigation. 
Such conserved quantities can be used to lower
the order of the Hamiltonian system of differential equations
and simplify the resolution of the optimal control problem.
They are also important for many other reasons.
In the calculus of variations they have been used to
synthesise the Lavrentiev  phenomenon \cite{MR90a:49020}
while in control, to analyze the stability and controllability
of nonlinear control systems, they are used for the system
decomposition in terms of simpler lower dimensional subsystems
\cite{MR86e:93030}. These are just few examples, but 
many other applications are possible: proving existence
of minimizers, solving the Hamilton-Jacobi-Bellman equation, etc.
We show that conserved quantities along the extremals 
are also a useful tool to prove Lipschitzian regularity
of the minimizing trajectories.

We will obtain some generalizations
of the well known theorems of E.~Noether, providing
a connection between such conserved quantities and the
invariance of the problems in optimal control. The theory
of this connection, as it appears in many branches of 
classical theoretical physics, constitutes one of the most
beautiful chapters of the calculus of variations.

The universal principle described by Noether's 
theorems of 1918, asserts that the invariance of a
problem with respect to a family of transformations
implies the existence of conserved quantities along
the Euler-Lagrange extremals. 
\begin{definition}
\label{Torres:defClassInv1}
If $C^1 \ni h^s(t,x) = \left(h_t^s(t),h_x^s(x)\right) : [a,b] \times \mathbb{R}^n 
\rightarrow \mathbb{R} \times \mathbb{R}^n$, 
$s \in \left(-\varepsilon,\varepsilon\right)$;
$h^0(t,x) = (t,x)$ for all $(t,x) \in [a,b] \times \mathbb{R}^n$;
\begin{equation*}
\int_{h_t^0(a)}^{h_t^0(b)}
L\left(t^s,h_x^s\left(x(t^s)\right),
\frac{\mathrm{d}}{\mathrm{d}t^s}h_x^s\left(x(t^s)\right)\right)\mathrm{d}t^s
= \int_a^b L\left(t,x(t),\dot{x}(t)\right)\mathrm{d}t \, ,
\end{equation*}
for $t^s = h_t^s(t)$, all $s \in (-\varepsilon,\varepsilon)$, and all $x(\cdot)$; 
then the basic problem of the calculus of variations is said to be invariant
under $h^s$.
\end{definition}
\begin{theorem}[First Noether's Theorem]
\label{Torres:ClassFNT}
If the basic problem of the calculus of variations is invariant under $h^s$, then
\begin{equation}
\label{Torres:eq:ClassFNT}
\psi(t) \cdot \frac{\partial}{\partial s}
\left.h_x^s\left(x(t)\right)\right|_{s = 0}
- H\left(t,x(t),\dot{x}(t),\psi(t)\right)
\frac{\partial}{\partial s}
\left.h_t^s\left(t\right)\right|_{s = 0}
\end{equation}
is constant in $t$ along every extremal.
\end{theorem}
We recall that for the basic problem of the
calculus of variations one has $\psi_0 = -1$ and 
$\psi = \frac{\partial L}{\partial u}$. Quantity
\eqref{Torres:eq:ClassFNT} is then equivalent to
\begin{multline}
\label{Torres:eq:equivClassFNT}
\frac{\partial L}{\partial u}\left(t,x(t),\dot{x}(t)\right) \cdot 
\frac{\partial}{\partial s}\left.h_x^s\left(x(t)\right)\right|_{s = 0} \\
+ \left[L\left(t,x(t),\dot{x}(t)\right) 
- \frac{\partial L}{\partial u}\left(t,x(t),\dot{x}(t)\right) \cdot \dot{x}(t)\right]
\frac{\partial}{\partial s} \left.h_t^s\left(t\right)\right|_{s = 0} \, .
\end{multline}
If the Lagrangian $L$
does not involve the time variable $t$ explicitly, 
one has invariance relative to translation with respect to time:
one can choose $h_t^s(t) = t + s$ and $h_x^s(x) = x$ in 
Definition~\ref{Torres:defClassInv1}.
It follows from Theorem~\ref{Torres:ClassFNT} that the corresponding Hamiltonian is
a first integral of the Euler-Lagrange equations. At the light of
\eqref{Torres:eq:equivClassFNT}, this is nothing more than the
classical $2^{nd}$ Erdmann condition,
\begin{equation}
\label{Torres:2ErdCond}
L\left(x(t),\dot{x}(t)\right) 
- \frac{\partial L}{\partial u}\left(x(t),\dot{x}(t)\right) \cdot \dot{x}(t)
\equiv \mbox{constant} \, ,
\end{equation}
which is a first-order
necessary optimality condition for the autonomous
basic problem of the calculus of variations.
Condition \eqref{Torres:2ErdCond} can also be obtained as a straight corollary
from the DuBois-Reymond necessary condition \eqref{Torres:eq:DuBoisReymond}
discussed in Section~\ref{Torres:Sec:CE}, and one can already guess the interplay
between the concept of Carath\'{e}odory equivalence and Noether theorems.
Such relation is the central key to obtain our generalizations. This is in
contrast with the classical proof of Theorem~\ref{Torres:ClassFNT}, 
which is based on the so called ``general variational formula''
\cite[pp. 172--198]{MR98b:49002a}. We claim that the proof 
based on Carath\'{e}odory approach is more far-reaching
than the traditional procedure. We will obtain a version
of Theorem~\ref{Torres:ClassFNT} to the optimal control setting
with several extensions and improvements 
(\textrm{cf.} Theorem~\ref{Torres:Cap5:r:mainresult} below).

Theorem~\ref{Torres:ClassFNT} comprises all theorems on first integrals
known to classical and quantum mechanics, field theories, and has deep implications 
in the general theory of relativity. For example, 
in mechanics \eqref{Torres:2ErdCond} correspond to the energy integral 
of conservative systems, a conservation law first discovered by Leonhard Euler in 1744;
while applying Noether's principle
to the Lagrangian describing a system of point masses,
one obtains conservation of linear momentum 
or angular momentum, corresponding, respectively, 
to invariance under spatial translation or spatial rotation.\footnote{In this
context the Hamiltonian multiplier $\psi(\cdot)$ represent the generalized
momentum of the system.}

As already mentioned, Noether's theorem can be considered
as an universal principle. Well known classical formulations
include invariant problems of the calculus of variations defined on a manifold
$M$; problems of the calculus of variations with multiple
integrals; invariance notions with respect to more than
one parameter; families of maps depending upon
arbitrary functions (second Noether theorem);
invariance of the Lagrangian up to addition of an
exact differential $\mathrm{d}\Phi(t,x,s)$, with
$\Phi$ linear on the parameter $s$.
In the original paper \cite{MR53:10538}, 
Noether explains that the derivatives 
of the state trajectories $x(\cdot)$ may also occur in the 
family of transformations $h^s$. However, this possibility
has been forgotten in the literature of the calculus of variations.
From our point of view, this possibility is very interesting:
it means that the parameter transformations $h^s$ may also depend
on the control variables.

Recent formulations, in other contexts than the calculus of variations, 
include the ones obtained by van der Schaft \cite{MR83k:49054} for 
autonomous Hamiltonian control systems with inputs and outputs; 
the results of Cari\~{n}ena and Figueroa \cite{MR96g:58011}
for (higher-order) supermechanics; and the discrete versions, 
in which time proceeds in integer steps, obtained
by Baez and Gilliam \cite{MR95i:58098}. 
In the optimal control setting, the important relation between 
invariance of the problem under 
a parameter family of transformations, and the existence 
of preserved quantities along the Pontryagin extremals, 
was established by Djukic \cite{MR49:5979}, 
Sussmann \cite{MR96i:49037}, Jurdjevic \cite[Ch. 13]{MR98a:93002}, 
Blankenstein \& van der Schaft \cite{MR1806135},
and Torres \cite{delfimEJC}.
Our purpose here is to provide 
a Noether type theorem to generic problems of optimal control
in a broader sense, enlarging the scope of its application.
Our results will be
formulated under a weak notion of invariance which admits
several parameters, equalities up to first order terms in the parameters, 
addition of an exact differential, not necessarily linear with
respect to the parameters, and a family of transformations which may 
also depend on the control variables. We will make use of a
technique different from the classical one.
This technique, introduced by the author in
\cite{delfimEJC}, does not need to use transversality conditions
as happens in the classical proof. For this reason,
the results will be valid even in the situation when 
we do not know the boundary conditions. 

We begin with a Noether theorem with no transformation
of the time variable.
\begin{definition}
\label{Torres:Cap5:d:invariancia}
Let $h^s : [a,b] \times \mathbb{R}^n \times \mathbb{R}^{r}
\rightarrow \mathbb{R}^n$,
\begin{equation*}
s = \left(s_1, \ldots, s_\rho\right) \, , 
\left\|s\right\| = \sqrt{\sum_{k=1}^{\rho} \left(s_k\right)^2} < \varepsilon \, ,
\end{equation*}
be a $\rho$-parametric family of $C^1$ transformations
which for $s = 0$ reduce to identity:
\begin{equation*}
h^0(t,x,u) = x \, , \quad \forall (t,x,u) \in [a,b] \times \mathbb{R}^n
\times \mathbb{R}^{r} \, .
\end{equation*}
If there exists a function $\Phi^s(t,x,u) \in
C^1\left([a,b],\mathbb{R}^n,\mathbb{R}^{r};\,\mathbb{R}\right)$ and
for all $s = \left(s_1, \ldots, s_\rho\right)$,  
$\left\|s\right\| < \varepsilon$, and admissible
$(x(\cdot),u(\cdot))$ there exists
a control $u^s(\cdot) \in L_{\infty}\left([a,b];\,\mathbb{R}^{r}\right)$
such that:
\begin{multline}
\label{Torres:Cap5:eq:inv1L}
\int_{a}^{\beta} L\left(t,h^s(t,x(t),u(t)),u^s(t)\right)\,\mathrm{d}t\\
= \int_{a}^{\beta} \left(L(t,x(t),u(t))
+ \frac{\mathrm{d}}{\mathrm{d} t}\Phi^s(t,x(t),u(t)) 
+ \delta\left(t,x(t),u(t),s\right) \right) \,\mathrm{d}t
\end{multline}
for all $\beta \in [a,b]$;
\begin{equation}
\label{Torres:Cap5:eq:inv1F}
\frac{\mathrm{d}}{\mathrm{d}t} h^s(t,x(t),u(t)) + \delta\left(t,x(t),u(t),s\right) =
\varphi(t,h^s(t,x(t),u(t)),u^s(t)) \, ;
\end{equation}
where $\delta\left(t,x,u,s\right)$ denote terms which go to zero faster than
$\left\|s\right\|$ for each $t$, $x$, $u$, \textrm{i.e.},
\begin{equation*}
\lim_{\left\|s\right\|\rightarrow 0} \frac{\delta\left(t,x,u,s\right)}{\left\|s\right\|} 
= 0 \, ,
\end{equation*}
then problem $(P)$ is said to be
quasi-invariant under $h^s(t,x,u)$ up to $\Phi^s(t,x,u)$.
\end{definition}
The following examples illustrate some of the new possibilities.

As far as in Definition~\ref{Torres:Cap5:d:invariancia}
$\Phi^s$ may depend on the parameters in a
nonlinear way, one may cover new situations
even for the basic problem of the calculus of variations.
Example~\ref{Torres:ex:pbcvTN} illustrate this issue.
We also note that in the example the state-transformation 
$h^s$ depend not only on the state variable
$x$ but also on time $t$. 
\begin{example}[$n = r = 1$]
\label{Torres:ex:pbcvTN}
Consider the following basic problem of the calculus
of variations:
\begin{gather*}
\int_a^b \left(u(t)\right)^2\,\mathrm{d}t \longrightarrow \min \, , \\
\dot{x}(t) = u(t)\, .
\end{gather*}
In this case we have $L = u^2$ and $\varphi = u$. The problem is invariant, 
in the sense of Definition~\ref{Torres:Cap5:d:invariancia},
under the one-parameter transformation $h^s(t,x) = x + s t$ ($h^0(t,x) = x$).
Indeed, we observe that for $\Phi^s(t,x) = s^2 t + 2 s x$ and
$u^s(t) = u(t) + s$ ($u^0(t) = u(t)$) one obtains:
\begin{gather*}
\begin{split}
\int_{a}^{\beta} L\left(u^s(t)\right)\, \mathrm{d}t
& = \int_{a}^{\beta} \left(u(t)+s\right)^2\, \mathrm{d}t
= \int_{a}^{\beta} \left(\left(u(t)\right)^2 + s^2 + 2 s u(t)\right)\, \mathrm{d}t \\
& = \int_{a}^{\beta} \left(L\left(u(t)\right) + 
\frac{\mathrm{d}}{\mathrm{d}t} \, \Phi^s\left(t,x(t)\right)\right)\, \mathrm{d}t \, ;
\end{split} \\
\frac{\mathrm{d}}{\mathrm{d}t} \, h^s\left(t,x(t)\right) = 
\dot{x}(t) + s = \varphi\left(u^s(t)\right)\, .
\end{gather*}
\end{example}

The optimal control problem in Example~\ref{Torres:ex:invFO}
is quasi-invariant under
a one-parameter family of transformations up
to an exact differential. 
\begin{example}[$n = 3$, $r = 2$] 
\label{Torres:ex:invFO}
Let us consider the following optimal control problem:
\begin{gather*} 
\int_a^b \left(u_1(t)\right)^2 + \left(u_2(t)\right)^2 \mathrm{d} t
\longrightarrow \min \, , \\
\begin{cases}
\dot{x_1}(t) = u_1(t) \, , \\
\dot{x_2}(t) = u_2(t)  \, , \\
\dot{x_3}(t) = \displaystyle \frac{u_2(t) \left(x_2(t)\right)^2}{2} \, .
\end{cases}
\end{gather*}
One has $L\left(u_1,u_2\right)=u_1^2 + u_2^2$,
$\varphi\left(x_2,u_1,u_2\right) = \left(u_1,u_2,\frac{u_2 x_2^2}{2}\right)^T$ and
direct computations show that the problem is
quasi-invariant under
\begin{equation*}
\begin{split}
h^s(t,x_1,x_2,x_3) &= \left(h_{x_1}^s\left(t,x_1\right),%
h_{x_2}^s\left(t,x_2\right),h_{x_3}^s\left(t,x_2,x_3\right)\right) \\
&= \left(x_1+st,x_2+st,x_3+\frac{1}{2} x_2^2 s t\right)
\end{split}
\end{equation*}
up to $\Phi^s(x_1,x_2) = 2s\left(x_1 + x_2\right)$:
$h^0(t,x_1,x_2,x_3) = (x_1,x_2,x_3)$ 
and by choosing $u_1^s = u_1+s$ and $u_2^s = u_2+s$ ($u_1^0 = u_1$, $u_2^0 = u_2$) 
we obtain
\begin{equation*}
\begin{split}
\int_a^\beta L\left(u_1^s(t),u_2^s(t)\right) \mathrm{d} t &= 
\int_a^\beta \left(u_1(t)+s\right)^2+\left(u_2(t)+s\right)^2 \mathrm{d} t \\
&= \int_a^\beta \left[\left(u_1(t)^2 + u_2(t)^2\right)  + 
2s\left(u_1(t) + u_2(t)\right) + 2 s^2 \right] \mathrm{d} t \\
&= \int_a^\beta \left[L\left(u_1(t),u_2(t)\right) 
+ \frac{\mathrm{d}}{\mathrm{d} t} \left(\Phi^s\left(x_1(t),x_2(t)\right)\right) 
+ \delta(s)\right] \mathrm{d} t \, ,
\end{split}
\end{equation*}
that is, equation \eqref{Torres:Cap5:eq:inv1L} is satisfied with $\delta(s) = 2 s^2$;
\begin{equation*}
\begin{split}
\varphi_1\left(u_1^s(t)\right) &= u_1(t) + s = \frac{\mathrm{d}}{\mathrm{d} t} (x_1(t) + st) 
= \frac{\mathrm{d}}{\mathrm{d} t} h_{x_1}^s\left(t,x_1(t)\right) \, , \\
\varphi_2\left(u_2^s(t)\right) &= u_2(t) + s = \frac{\mathrm{d}}{\mathrm{d} t} (x_2(t) + st) 
= \frac{\mathrm{d}}{\mathrm{d} t} h_{x_2}^s\left(t,x_2(t)\right) \, ,
\end{split}
\end{equation*}
\begin{equation*}
\begin{split}
\varphi_3 &\left(h_{x_2}^s\left(t,x_2(t)\right),u_2^s(t)\right) = 
\frac{\left(u_2(t) + s\right)\left(x_2(t)+st\right)^2}{2} \\
&= \frac{u_2(t) x_2(t)^2}{2} + \frac{1}{2} s \left(x_2(t)^2 + 2 x_2(t) u_2(t) t\right)
+ \frac{\left(u_2(t) t^2 + 2 x_2(t) t\right) s^2 + t^2 s^3}{2} \\
&= \frac{\mathrm{d}}{\mathrm{d} t}\left(x_3(t)+\frac{1}{2} x_2(t)^2 s t\right) 
+ \delta\left(t,x_2(t),u_2(t),s\right) \\
&= \frac{\mathrm{d}}{\mathrm{d} t} h_{x_3}^s\left(t,x_2(t),x_3(t)\right) 
+ \delta\left(t,x_2(t),u_2(t),s\right) \, ,
\end{split}
\end{equation*}
and therefore \eqref{Torres:Cap5:eq:inv1F} is also satisfied.
\end{example}

\begin{theorem}
\label{Torres:Cap5:r:cpten}
If $(P)$ is quasi-invariant
under the transformations $h^s(t,x,u)$ up to
$\Phi^s(t,x,u)$, in the sense of Definition~\ref{Torres:Cap5:d:invariancia}, then
the $\rho$ quantities
\begin{equation*}
\psi(t) \cdot \frac{\partial}{\partial s_k}
\left.h^s(t,x(t),u(t))\right|_{s=0} + \psi_0 \frac{\partial}{\partial s_k}
\left.\Phi^s(t,x(t),u(t))\right|_{s = 0} \, , \quad (k = 1,\ldots,\rho) \, ,
\end{equation*}
are constant in $t$ along every Pontryagin extremal
$\left(x(\cdot),u(\cdot),\psi_0,\psi(\cdot)\right)$ of the problem.
\end{theorem}

\begin{example}
From Theorem~\ref{Torres:Cap5:r:cpten} it follows that 
$\psi(t) t + 2 \psi_0 x(t)$ is
constant in $t$ along the extremals of the basic problem
of the calculus of variations considered in Example~\ref{Torres:ex:pbcvTN}.
We know that 
$\psi(t) = \frac{\partial L}{\partial u}\left(t,x(t),u(t)\right)$ 
and $\psi_0 = -1$, so the conclusion is that $\dot{x}(t) t - x(t)$ 
is a first integral of the Euler-Lagrange differential equations.
\end{example}

\begin{example}
One concludes from Theorem~\ref{Torres:Cap5:r:cpten} that
$2 \psi_0\left(x_1(t) + x_2(t)\right) + \psi_1(t) t + \psi_2(t) t +
\frac{1}{2} \psi_3(t) \left(x_2(t)\right)^2 t$ is constant
along all the Pontryagin extremals of the problem considered in 
Example~\ref{Torres:ex:invFO}.
\end{example}

\begin{example}
Problem \eqref{Torres:LuFiuvezx} is invariant under
the one-parameter family of transformations 
$h^s(x) = \mathrm{e}^s x$ (\textrm{cf.} 
Definition~\ref{Torres:Cap5:d:invariancia}
with $u^s = u$ and $\Phi^s \equiv 0$). We get from 
Theorem~\ref{Torres:Cap5:r:cpten} that 
\begin{equation}
\label{Torres:eq:LC1}
\psi(t) x(t) \equiv \mbox{constant} \, ,
\end{equation}
$t \in [a,b]$, along any Pontryagin extremal of problem
\eqref{Torres:LuFiuvezx}.
\end{example}

\begin{example}
\label{Torres:ex:34}
Under the condition \eqref{Torres:eq:homogeneidade}
problem \eqref{Torres:eq:pr:cubicPolyRiemMani} is invariant
under $h_{x_1}^s = \mathrm{e}^s x_1$, $h_{x_2}^s = \mathrm{e}^s x_2$,
with $u_i^s = u_i$ and $\Phi^s \equiv 0$. 
Proposition~\ref{Torres:prop:homogenei} 
follows from Theorem~\ref{Torres:Cap5:r:cpten}. Other conclusions
are also possible if one imposes different conditions on the vector fields 
$X_i$ (\textrm{cf.} Example~\ref{Torres:ex:MgCam}).
\end{example}

We now generalize the invariance notion
given by Definition~\ref{Torres:Cap5:d:invariancia},
in order to admit the possibility of a
$\rho$-parametric transformation
of the independent variable $t$.

\begin{definition}
\label{Torres:Cap5:d:invarianciaTX}
Let $h^s(t,x,u)=\left(h_t^s(t,x,u),\,h_x^s(t,x,u)\right)$,
$h_t^s : [a,b] \times \mathbb{R}^n \times \mathbb{R}^r
\rightarrow \mathbb{R}$ and
$h_x^s : [a,b] \times \mathbb{R}^n \times \mathbb{R}^r
\rightarrow \mathbb{R}^n$
($\left\|s\right\| < \varepsilon$),
be a $\rho$-parametric family of $C^1$ transformations which
for $s = 0$ satisfies $h^0(t,x,u) = (t,x)$ for all triple $(t,x,u) \in
[a,b] \times \mathbb{R}^n \times \mathbb{R}^r$.
If there exists a function $\Phi^s(t,x,u) \in
C^1\left([a,b],\mathbb{R}^n,\mathbb{R}^r;\,\mathbb{R}\right)$ and
for all $s = \left(s_1,\ldots,s_\rho\right)$ and for all
admissible pair $\left(x(\cdot),u(\cdot)\right)$ there exists a
control $u^s(\cdot) \in L_{\infty}\left([a,b];\,\mathbb{R}^r\right)$
such that:
\begin{multline*}
\int_{h_t^s(a,x(a),u(a))}^{h_t^s(\beta,x(\beta),u(\beta))} L\left(t^s,
h_x^s(t^s,x(t^s),u(t^s)),u^s(t^s)\right)\mathrm{d} t^s \\
= \int_{a}^{\beta} \left( L(t,x(t),u(t))
+ \frac{\mathrm{d}}{\mathrm{d} t}\Phi^s(t,x(t),u(t)) 
+ \delta\left(t,x(t),u(t),s\right)\right)\mathrm{d}t\, ,
\end{multline*}
for all $\beta \in [a,b]$;
\begin{equation*}
\frac{\mathrm{d}}{\mathrm{d} t^s} h_x^s(t^s,x(t^s),u(t^s)) 
+ \delta\left(t,x(t),u(t),s\right)
= \varphi(t^s,\,h_x^s(t^s,x(t^s),u(t^s)),\,u^s(t^s))\, ,
\end{equation*}
for $t^s = h_t^s(t,x(t),u(t))$;
then the problem $(P)$ is said to be quasi-invariant under
transformations $\left(h_t^s(t,x,u),\,h_x^s(t,x,u)\right)$
up to $\Phi^s(t,x,u)$.
\end{definition}

What follows is a more general version of the first
Noether theorem which admits transformations 
of the time-variable.
Theorem~\ref{Torres:Cap5:r:mainresult} 
gives $\rho$ conservation laws when the
optimal control problem $(P)$ is quasi-invariant
up to $\Phi^s$ under a family
of transformations with $\rho$ parameters.

\begin{theorem}[First Noether Theorem for Optimal Control]
\label{Torres:Cap5:r:mainresult}
If $(P)$ is quasi-invariant under the transformations
$\left(h_t^s(t,x,u),h_x^s(t,x,u)\right)$
up to $\Phi^s(t,x,u)$,
in the sense of Definition~\ref{Torres:Cap5:d:invarianciaTX}, then
\begin{multline*}
\psi(t) \cdot \frac{\partial}{\partial s_k}
\left.h_x^s(t,\,x(t),\,u(t))\right|_{s = 0}
+ \psi_0 \frac{\partial}{\partial s_k}
\left.\Phi^s(t,\,x(t),\,u(t))\right|_{s = 0} \\
- H\left(t,\,x(t),\,u(t),\,\psi_0,\,\psi(t)\right)
\frac{\partial}{\partial s_k} \left.h_t^s(t,\,x(t),\,u(t))\right|_{s = 0}
\end{multline*}
is constant in $t$ along every Pontryagin extremal
$\left(x(\cdot),u(\cdot),\psi_0,\psi(\cdot)\right)$ of the problem
$\forall \, k=1,\ldots,\rho$.
\end{theorem}

The proof of Theorem~\ref{Torres:Cap5:r:mainresult} is done by reduction to
the situation of Theorem~\ref{Torres:Cap5:r:cpten}. 
If $(P)$ is quasi-invariant under
$\left(h_t^s(t,x,u),h_x^s(t,x,u)\right)$
up to $\Phi^s(t,x,u)$, in the sense
of Definition~\ref{Torres:Cap5:d:invarianciaTX}, then
the problem $\left(P_\tau\right)$ introduced in Section~\ref{Torres:Sec:CE}
is quasi-invariant under $h^s\left(t,z,w\right)=
\left(h_t^s(t,z,w),h_x^s(t,z,w)\right)$ 
up to $\Phi^s(t,z,w)$
in the sense of Definition~\ref{Torres:Cap5:d:invariancia}.
The pretended conclusion is then obtained from 
Theorem~\ref{Torres:Cap4:r:prop2}.

\begin{remark}
Every autonomous problem is invariant under
$h_t^s = t + s$ and $h_x^s = x$ 
(autonomous problems are time-invariant). 
It follows from Theorem~\ref{Torres:Cap5:r:mainresult}
that the corresponding Hamiltonian $H$ is
constant along the Pontryagin extremals
(\textrm{cf.} equality \eqref{Torres:eq:dHdtpartialHdt}).
For the time-invariant problem \eqref{Torres:LuFiuvezx},
this fact, together with the conservation law \eqref{Torres:eq:LC1},
gives an alternative explanation for \eqref{Torres:eq:CLCh1}
to be constant along the extremals of the problem.
\end{remark}

\begin{example}
\label{Torres:ex:MgCam}
Under the hypotheses 
$X_i\left(\lambda x_1\right) = \lambda^{\alpha} X_i\left(x_1\right)$,
$\alpha \in \mathbb{R}\setminus \left\{1\right\}$,
problem \eqref{Torres:eq:pr:cubicPolyRiemMani} is invariant under 
\begin{gather*}
t^s = h_{t}^s(t) = \mathrm{e}^{-2s} t \, , \,
h_{x_1}^s\left(x_1(t^s)\right) = \mathrm{e}^{\frac{3}{\alpha - 1}s} x_1(t) \, , \,
h_{x_2}^s\left(x_2(t^s)\right) = 
\mathrm{e}^{\left(\frac{3\alpha}{\alpha - 1} - 1\right)s} x_2(t) \, ,
\end{gather*}
with $u_i^s(t^s) = \mathrm{e}^s u_i(t)$ and $\Phi^s \equiv 0$. We conclude from 
Theorem~\ref{Torres:Cap5:r:mainresult} that
\begin{equation*}
\psi_1(t) \frac{3}{\alpha - 1} x_1(t) 
+ \psi_2(t) \left(\frac{3\alpha}{\alpha - 1} - 1\right) x_2(t)
+ 2 H t \equiv \mbox{constant}
\end{equation*}
holds along any extremal of the problem \eqref{Torres:eq:pr:cubicPolyRiemMani}.
\end{example}

It is interesting to note that Theorem~\ref{Torres:Cap5:r:mainresult} 
cover both normal and abnormal situations. This is an important
issue because abnormal minimizers while nonexistent
for the basic problem of the calculus of variations, 
in general Lagrange problem they may
occur frequently. This is the case, for example,
for the problems in Sub-Riemannian Geometry.

Even for the basic problem of the calculus of variations
the results are new and provide new information.
As far as the notions of invariance,
conserved quantity along the extremals and reduction belong to the
most important tools in the study of classical mechanics,
it is not surprising that the conserved quantities obtained
by Theorem~\ref{Torres:Cap5:r:mainresult} may be very useful in practice.
We remark that solving the Hamiltonian system by 
the elimination of the control, with the aid of the maximality condition,
is typically a difficult task. The existence 
of such conserved quantities are a
circumstance which may make the resolution process easier
and are often useful for purposes of analyzing a nonlinear control system.
We shall see in the next section that the conserved quantities obtained
by Theorem~\ref{Torres:Cap5:r:mainresult} are also useful
to establish Lipschitzian regularity of the minimizing trajectories.

%%%%%%%%%%%%%%%%%%%%%%%%%%%

\section{Tonelli Full-Regularity}
\label{Torres:Sec:TFR}

In this section we address the question: 
\emph{How to assure that the set of extremals include the minimizers predicted 
by the existence theory?}

It is easy to find examples of the optimal control problem $(P)$, 
very simple in aspect, for which the application of the
Pontryagin maximum principle gives a unique function $u(\cdot)$
which is not an optimal control. This happens
because the optimal solution does not exist. One cannot
conclude that we have found the solution unless we know 
\emph{a priori} that a solution really exists. 

A general existence theory for the calculus
of variations has been introduced by Leonida Tonelli,
in a series of Italian papers, as from 1911, when he was 26.
The first general existence theorem for optimal control
was given by Filippov. The original paper, in Russian, 
appeared in 1959. 
There exist now an extensive literature on the existence of
solutions to problems of optimal control.
The following set of conditions, of the type of
Tonelli, guarantee the existence of minimizer for our problem $(P)$.
It is called a Tonelli type existence theorem because
for the basic problem of the calculus of variations
one has $\varphi = u$ and the theorem
coincides with the classical Tonelli existence theorem.
\begin{theorem}[``Tonelli'' Existence Theorem for $(P)$]
Problem $(P)$ has an absolute minimum in the space
$u(\cdot) \in L_1$, provided that there exist at least one admissible pair, 
and the following conditions are satisfied for all $\left(t,x,u\right)$:
\begin{itemize}
\item Coercivity: there exists a function
$\theta : \mathbb{R}_{0}^{+}\rightarrow \mathbb{R}$, bounded below, such that
\begin{gather*}
\lim_{r\rightarrow +\infty}\frac{\theta(r)}{r} = + \infty \, , \\
L\left(t,x,u\right) \geq \theta\left(\left\|\varphi(t,x,u)\right\|\right) \, , \\
\lim_{\left\|u\right\| \rightarrow + \infty} 
\left\|\varphi(t,x,u)\right\| = + \infty \, ;
\end{gather*}
\item Convexity: $L\left(t,x,u\right)$ and
$\varphi\left(t,x,u\right)$ are convex with respect to $u$.
\end{itemize}
\end{theorem}
Roughly speaking, the theorem asserts that under convexity and coercivity,
a solution exists in the class of integrable controls.
For the basic problem of the calculus of
variations one has $\dot{x} =u$ and
this mean that existence is
given in the class of absolutely continuous
functions, possible with unbounded derivative.
We note that the assumptions on the solution
for the derivation of the necessary optimality conditions
have more regularity than the one considered here.
For example, for the basic problem of the calculus of
variations, the biggest class for which the Euler-Lagrange
equation is valid is the class of Lipschitzian functions,
that is, the class of absolutely continuous functions having
essentially bounded derivative.
The steps for the derivation of the the Euler-Lagrange equations 
can no longer be justified in the class of absolutely continuous functions.
So the central question, which immediately comes to mind, is the following:
How different is the problem with controls in $L_1$
from the problem with the controls in $L_\infty$?
It seems that in order to apply the standard approach 
to solving optimization problems,
one needs an intermediate step between existence,
which is proved for controls in $L_1$,
and standard classical necessary conditions, which are valid
for controls in $L_\infty$.
Is this intermediate step really necessary?
Is this a technical phenomenon or does it reflect a fundamental
difficulty? Can the solution predicted by Tonelli's existence theorem
be irregular and fail to satisfy the Pontryagin maximum principle?
Tonelli proved, for the basic problem of the calculus of variations 
in the scalar case ($n = 1$), that bad behaviour is only possible
in a closed set of measure zero.
It turns out that, as has been shown by F.~H.~Clarke and R.~B.~Vinter
\cite[Ch. 11]{MR2001c:49001},
that the result is general. 
\begin{theorem}[``Tonelli'' Regularity]
Assume that the Tonelli Existence Hypotheses are satisfied.
Take any minimizer $\left(\tilde{x}(\cdot),\tilde{u}(\cdot)\right)$
of $(P)$. Then there exists a closed subset $\Omega \subset [a,b]$
of zero measure with the following property: for any
$\tau \in [a,b] \setminus \Omega$, $\tilde{u}(\tau)$ is
essentially bounded on a relative neighborhood of $\tau$.
\end{theorem}
The theorem asserts that the optimal solutions predicted by Tonelli's
existence theorem satisfy the Pontryagin maximum principle
or the Euler-Lagrange equations everywhere, except possibly
at the points of a closed exceptional set $\Omega$ of measure zero.
Tonelli, against the general opinion that, 
at least for ``reasonable'' problems,
all minimizers predicted by his existence
theorem are extremals, conjectured the possibility
of the $\Omega$ set to be nonempty and the possibility
of such a minimizer $\tilde{u}(\cdot)$ to be unbounded
and fail to be an extremal.
It has been proved in 1984 by J. M. Ball and V. J. Mizel 
\cite{MR86k:49002} that
in general the $\Omega$ set can not be taken to be empty,
even for very ``reasonable'' problems. Even for polynomial Lagrangians and 
linear dynamics, minimizers predicted by Tonelli's existence theorem
may fail to be Pontryagin (or Euler-Lagrange) extremals.
Given this possibility, the natural question to ask now is
the following: How to exclude the possibility of bad behaviour?
How to obtain full-regularity ($\Omega = \emptyset$)? 
This is achieved by postulating conditions beyond those of 
Tonelli's existence theorem, assuring that
all optimal controls are essentially bounded.
These conditions close the gap between the hypotheses
arising in the existence theory and those of necessary optimality
conditions, assuring that the solutions can be identified via the 
Pontryagin maximum principle. As far as $\varphi$ is bounded on bounded
sets, it also follows that the optimal trajectory is
Lipschitzian and, similarly, the Hamiltonian adjoint multipliers
$\psi$ of the Pontryagin maximum principle turn out to
be Lipschitzian either. Thus, full-regularity justifies
searching for minimizers among extremals and establishes a weaker
form of the maximum principle in which the Hamiltonian adjoint 
multipliers are not required to be absolutely continuous but
merely Lipschitzian. With the Lipschitzian regularity in hand, 
other regularity properties follow easily, like $C^1$ or $C^2$,
or even $C^\infty$, imposing some more additional conditions. 
Full-regularity conditions also precludes occurrence of
undesirable phenomena, like the Lavrentiev one,
making possible the implementation of efficient discretization
schemes and algorithms for numerical
computation of the optimal controls. Again,
such undesirable phenomena are
possible even when the Lagrangian is a polynomial
and the control system is linear \cite{mania}.

The regularity condition one most often finds,
implying that the minimizing controls are bounded, 
was suggested by Tonelli himself.
Tonelli and Morrey proved that under the growth conditions
\begin{equation}
\label{Torres:eq:TMTC}
\left\| L_{x}\right\| +\left\| L_{u}\right\| \leq c\,\left| L\right| + k \, ,
\end{equation}
$c>0$ and $k$ constants, 
points of bad behaviour cannot occur for the
basic problem of the calculus of variations:
under conditions \eqref{Torres:eq:TMTC} all optimal controls predicted
by Tonelli's existence theorem are bounded and
the corresponding minimizing trajectories Lipschitzian.
The conditions impose global growth hypothesis
on the derivatives of the Lagrangian $L$ with respect
to the state and control variables. 
F.~H.~Clarke and R.~B.~Vinter have shown that the bound on the
derivatives of the Lagrangian with respect to the control
variables can be discarded, and that regularity conditions
\begin{equation*}
\left\|L_{x}\right\| \leq c \left|L\right| + k \, ,
\end{equation*}
of the type of Tonelli-Morrey, hold not only for the basic problem but
universally in the calculus of variations \cite[Ch. 11]{MR2001c:49001}.
We will see that Tonelli-Morrey-type regularity conditions apply in fact
more generally: they hold in the generic context of optimal control. 
They hold even when the dynamics are nonlinear both in the state 
and control variables.

The literature on regularity conditions for the problems
of the calculus of variations is now vast, but for the
problems of optimal control, if one excludes the special cases
that can be easily recast as problems in the calculus of
variations, regularity conditions are a rarity. 
The first results appeared in 2000 and treat the case of 
control-affine dynamics:
\begin{theorem}[\cite{MR2000m:49048}]
\label{Torres:th:MR2000m:49048}
For the Lagrange Problem of Optimal Control $(P)$ with
control affine dynamics, $\varphi = f(t,x) + g(t,x)\,u$,
if $g\left(t,x\right)$ has complete rank $r$
for all $t$ and $x$; the coercivity condition holds; and
$\exists$ $\gamma >0$, $\beta <2$, 
$\eta $, and $\mu \geq \max \left\{\beta -2,\,-2\right\}$, such that 
\begin{equation}
\label{Torres:RC:ASDT}
\left( \left| L_{t}\right| +\left| L_{x^{i}}\right| +\left\| L\,
\varphi_{t}-L_{t}\,\varphi \right\| +\left\| L\,\varphi _{x^{i}}-L_{x^{i}}\,
\varphi\right\| \,\right) \,\left\| u\right\| ^{\mu }\leq \gamma \,L^{\beta }
+\eta \, ,
\end{equation}
then all the minimizers $\tilde{u}\left(\cdot\right)$ of the problem, 
which are not abnormal extremal controls, are essentially bounded on 
$\left[a,b\right]$.
\end{theorem}
The proof of Theorem~\ref{Torres:th:MR2000m:49048} is based on the reduction of
problem $(P)$ to problem \eqref{Torres:eq:MudVarIndCG} with
$\Upsilon(t,z,v) = \frac{1}{L\left(t,z,v\right)}$;
on the subsequent Gamkrelidze's compactification of the space of 
admissible controls \cite[\S 8.5]{MR58:33350c}; 
on the abnormal-Carath\'{e}odory-equivalence
given by Theorems~\ref{Torres:Cap4:T1NP} and \ref{Torres:Cap4:ZeroLevel:OWA};
and utilization of the classic Pontryagin maximum principle and
the time-invariance property of problem \eqref{Torres:eq:MudVarIndCG}.
We remark that conditions \eqref{Torres:RC:ASDT} are not of the type of Tonelli-Morrey,
and even for the basic problem of the calculus of variations
one can cover new situations \cite{MR2000m:49048,MR2001j:49062}.
The only drawback is that in order to cover new situations the regularity conditions 
become harder to verify. For the case of control-affine dynamics this is not a problem,
and it is possible to deal pretty well with the conditions. 
In \cite{MR2000m:49048} other conditions, not so general as 
\eqref{Torres:RC:ASDT}, more strong, but more easy to check in practice,
were obtained. Theorem~\ref{Torres:th:MR2000m:49048} admit a generalization for Lagrange
problems with dynamics which are nonlinear in control, introducing
generalized controls and making a reduction of the nonlinear
dynamics to the control affine case by relaxation, a technique
introduced by R.~V.~Gamkrelidze. The only problem, with this nice
approach, is that the conditions become cumbersome.
We must not forget that checking regularity
conditions is a preliminary step in the process of solving
a problem, and that, by definition,
regularity conditions must be simple to verify.
To go to the general nonlinear case
with verifiable regularity conditions, a new technique is needed. 
Such technique was introduced by the author in \cite{torresMCSS},
showing that Tonelli-Morrey type conditions
work universally in optimal control:
\begin{theorem}[\cite{torresMCSS}]
\label{Torres:MR:MCSS}
Coercivity plus the growth conditions: there exist
constants $c>0$ and $k$ such that
\begin{equation}
\label{Torres:LRCTMT}
\begin{split}
\left|\frac{\partial L}{\partial t}\right| &\le
c \left|L\right| + k \, , \quad 
\left\|\frac{\partial L}{\partial x}\right\| \le
c \left|L\right| + k \, ,\\
\left\|\frac{\partial \varphi}{\partial t}\right\| &\le
c \left\|\varphi\right\| + k \, , \quad 
\left\|\frac{\partial \varphi_i}{\partial x}\right\| \le
c \left|\varphi_i\right| + k \quad (i = 1,\,\ldots,\,n) \, ;
\end{split}
\end{equation}
imply that all minimizers $\tilde{u}(\cdot)$ of $(P)$,
which are not abnormal extremal controls,
are essentially bounded on $[a,b]$.
\end{theorem}
The Lipschitzian regularity conditions \eqref{Torres:LRCTMT}
are obtained using the applicability conditions 
\eqref{Torres:CondApplPMPBerkClar} of the Pontryagin maximum principle 
to the auxiliary problem $(P_\tau)$ introduced in Section~\ref{Torres:Sec:CE}.
Theorem~\ref{Torres:MR:MCSS} is then proved using the theorem of
Emmy Noether and the established abnormal-Carath\'{e}odory-equivalence
between problems $(P)$ and $(P_\tau)$.
The theorem covers the general optimal control problem $(P)$,
providing conditions of the type of Tonelli-Morrey under
which non-abnormal optimal controls are bounded. This
guarantees the Lipschitzian regularity of the non-abnormal minimizing
trajectories and that all minimizers are Pontryagin extremals.
\begin{corollary}
Under the hypotheses of Theorem~\ref{Torres:MR:MCSS},
all minimizers of $(P)$ are Pontryagin extremals.
\end{corollary}
We remark that convexity is not required in 
Theorems~\ref{Torres:th:MR2000m:49048} and \ref{Torres:MR:MCSS} 
in order to establish the Lipschitzian regularity of the (non-abnormal)
minimizing trajectories $\tilde{x}(\cdot)$. This fact is
important because existence theorems without the convexity 
assumptions are a question of great interest.

It is also possible to obtain new regularity conditions, 
which are not of the type of Tonelli-Morrey, for the generic
nonlinear problem $(P)$ \cite{torresMCSS}.

I found pertinent to quote 
Constantin Carath\'{e}odory addressing the question of
\emph{The Beginning of Research in the Calculus
of Variations} \cite{ZBL0018.19601}: 
``I will be glad if I have succeeded in impressing the idea that 
it is not only pleasant and entertaining to read at times the works 
of the old mathematical authors, but that this may occasionally 
be of use for the actual advancement of science.''

%%%%%%%%%%%%%%%%%%%%%%%%%%%%%%%%%%%%%%%%%%%%%%%%%%%%%%%%%%%%%%%%%%%

\section*{Acknowledgments}

I would like to thank the support 
from Project POCTI/MAT/41683/2001
\emph{Advances in Nonlinear Control and Calculus of Variations},
FCT -- Sapiens'01, of the Research and Development Unit 
CEOC (Centro de Estudos em Optimiza\c{c}\~{a}o e Controlo),
University of Aveiro, Portugal.

%%%%%%%%%%%%%%%%%%%%%%%%%%%%%%%%%%%%%%%%%%%%%%%%%%%%%%%%%%%%%%%%%%%

%%%%%%%%%%%%%%%%%%%%%%%%%%%%%%%%%%%%%%%%%%%%%%%%%%%%%%%%%%%%%%%%%%%


\begin{thebibliography}{10}

\bibitem{MR95i:58098}
J.~C. Baez and J.~W. Gilliam.
An algebraic approach to discrete mechanics.
{\em Lett. Math. Phys.}, 31(3):205--212, 1994. 
\textsf{Zbl} {0805.58031} \textsf{MR} {95i:58098}

\bibitem{MR86k:49002}
J.~M. Ball and V.~J. Mizel.
One-dimensional variational problems whose minimizers do not satisfy
the {E}uler {L}agrange equation. 
{\em Arch. Rational Mech. Anal.}, 90(4):325--388, 1985.
\textsf{Zbl} {0585.49002} \textsf{MR} {86k:49002}

\bibitem{MR51:8914}
L.~D. Berkovitz.
{\em Optimal control theory}.
Springer-Verlag, New York, 1974. 
\textsf{Zbl} {0295.49001} \textsf{MR} {51:8914}

\bibitem{MR1806135}
G.~Blankenstein and A.~van~der Schaft.
Optimal control and implicit {H}amiltonian systems.
In {\em Nonlinear control in the year 2000, Vol.\ 1 (Paris)}, 
pages 185--205. Springer, London, 2001. 
\textsf{MR} {1806135}

\bibitem{ZBL0018.19601}
C.~Carath\'{e}odory.
The beginning of research in the calculus of variations.
{\em Osiris}, 3:224--240, 1937. 
\textsf{Zbl} {0018.19601}

\bibitem{ZBL0505.49001}
C.~Carath\'{e}odory.
{\em Calculus of variations and partial differential equations of the
first order}. Chelsea Publishing Company, New York, 1982. 
\textsf{Zbl} {0505.49001}

\bibitem{MR96g:58011}
J.~F. Cari\~nena and H.~Figueroa.
A geometrical version of {N}oether's theorem in supermechanics.
{\em Rep. Math. Phys.}, 34(3):277--303, 1994. 
\textsf{Zbl} {0846.58008} \textsf{MR} {96g:58011}

\bibitem{MR54:3540}
F.~H. Clarke.
The maximum principle under minimal hypotheses.
{\em SIAM J. Control Optimization}, 14(6):1078--1091, 1976.
\textsf{Zbl} {0344.49009} \textsf{MR} {54:3540}

\bibitem{MR1796845}
F.~H. Clarke.
The calculus of variations, nonsmooth analysis and optimal control.
In {\em Development of mathematics 1950--2000}, pages 313--328.
Birkh\"auser, Basel, 2000. 
\textsf{Zbl} {0970.49002} \textsf{MR} {2001h:49003}

\bibitem{MR58:33350c}
R.~V. Gamkrelidze.
{\em Principles of optimal control theory}.
Plenum Press, New York, 1978. 
\textsf{Zbl} {0401.49001} \textsf{MR} {58:33350c}

\bibitem{MR98b:49002a}
M.~Giaquinta and S.~Hildebrandt.
{\em Calculus of variations {I}. The Lagrangian formalism}.
Springer-Verlag, Berlin, 1996. 
\textsf{Zbl} {0853.49001} \textsf{MR} {98b:49002a}

\bibitem{MR86e:93030}
J.~W. Grizzle and S.~I. Marcus.
The structure of nonlinear control systems possessing symmetries.
{\em IEEE Trans. Automat. Control}, 30(3):248--258, 1985.
\textsf{MR} {86e:93030}

\bibitem{MR90a:49020}
A.~C. Heinricher and V.~J. Mizel.
The {L}avrentiev phenomenon for invariant variational problems.
{\em Arch. Rational Mech. Anal.}, 102(1):57--93, 1988.
\textsf{Zbl} {0655.49003} \textsf{MR} {90a:49020}

\bibitem{MR98a:93002}
V.~Jurdjevic.
{\em Geometric control theory}.
Cambridge University Press, Cambridge, 1997. 
\textsf{Zbl} {0940.93005} \textsf{MR} {98a:93002}

\bibitem{mania}
B.~Mani\`a.
Sopra un esempio di {L}avrentieff.
{\em Boll. Un. Mat. Ital.}, 13:147--153, 1934.

\bibitem{JFM46.0770.01}
E.~Noether.
Invariante variationsprobleme.
{\em G\"{o}tt. Nachr.}, pages 235--257, 1918. 
\textsf{JFM} {46.0770.01}

\bibitem{MR53:10538}
E.~Noether.
Invariant variation problems.
{\em Transport Theory Statist. Phys.}, 1(3):186--207, 1971.
English translation of the original paper \cite{JFM46.0770.01}.
\textsf{Zbl} {0292.49008} \textsf{MR} {53:10538}

\bibitem{MR29:3316b}
L.~S. Pontryagin, V.~G. Boltyanskii, R.~V. Gamkrelidze, and E.~F. Mishchenko.
{\em The mathematical theory of optimal processes}.
Interscience Publishers John Wiley \& Sons, Inc.\, New York-London, 1962. 
\textsf{Zbl} {0882.01027} \textsf{MR} {29:3316b}

\bibitem{MR2000m:49048}
A.~V. Sarychev and D.~F.~M. Torres.
Lipschitzian regularity of minimizers for optimal control problems
with control-affine dynamics.
{\em Appl. Math. Optim.}, 41(2):237--254, 2000. 
\textsf{Zbl} {0961.49021} \textsf{MR} {2000m:49048}

\bibitem{MR2001j:49062}
A.~V. Sarychev and D.~F.~M. Torres.
Lipschitzian regularity conditions for the minimizing trajectories of
optimal control problems.
In {\em Nonlinear analysis and its applications to differential
equations (Lisbon, 1998)}, pages 357--368. Birkh\"auser Boston, Boston, MA, 2001. 
\textsf{Zbl} {pre01693620} \textsf{MR} {2001j:49062}

\bibitem{MR96i:49037}
H.~J. Sussmann.
Symmetries and integrals of motion in optimal control.
In {\em Geometry in nonlinear control and differential inclusions
(Warsaw, 1993)}, pages 379--393. Polish Acad. Sci., Warsaw, 1995.
\textsf{Zbl} {0891.49011} \textsf{MR} {96i:49037}

\bibitem{MR1806192}
H.~J. Sussmann.
New theories of set-valued differentials and new versions of the
maximum principle of optimal control theory.
In {\em Nonlinear control in the year 2000, Vol.\ 2 (Paris)}, 
pages 487--526. Springer, London, 2001. 
\textsf{Zbl} {pre01584893} \textsf{MR} {2002e:49040}

\bibitem{300Sussmann}
H.~J. Sussmann and J.~C. Willems.
300 years of optimal control: from the brachystochrone to the maximum
principle. {\em IEEE Control Systems}, pages 32--44, 1997.

\bibitem{delfim3ncnw}
D.~F.~M. Torres.
Conservation laws in optimal control.
In {\em Dynamics, Bifurcations and Control}, volume 273 of 
{\em Lecture Notes in Control and Information Sciences}, pages 287--296.
Springer-Verlag, Berlin, Heidelberg, 2002.

\bibitem{torresMED2002}
D.~F.~M. Torres.
Conserved quantities along the {P}ontryagin extremals of
quasi-invariant optimal control problems.
In {\em Proc. 10th Mediterranean Conference on Control and
Automation, MED2002, Lisboa, Portugal}, 2002.

\bibitem{torresCM02I06E}
D.~F.~M. Torres.
On optimal control problems which admit an infinite continuous group
of transformations. In {\em Proc. 5th Portuguese Conference on Automatic Control,
Controlo 2002, Aveiro, Portugal}, 2002.

\bibitem{delfimEJC}
D.~F.~M. Torres.
On the {N}oether theorem for optimal control.
{\em European Journal of Control}, 8(1):56--63, 2002.

\bibitem{torresPhD}
D.~F.~M. Torres.
{\em Regularity of Minimizers in the Calculus of Variations and
Optimal Control}. {P}h.{D}. thesis, Dep. Mathematics, 
Univ. Aveiro, Portugal, 2002. (In Portuguese).

\bibitem{delfimIO}
D.~F.~M. Torres.
A remarkable property of the dynamic optimization extremals.
{\em Investiga\c{c}\~{a}o Operacional}, in press.

\bibitem{torresMCSS}
D.~F.~M. Torres.
Lipschitzian regularity of the minimizing trajectories for nonlinear
optimal control problems. Submitted for publication.

\bibitem{MR49:5979}
{\leavevmode\lower.6ex\hbox to 0pt{\hskip-.23ex \accent"16\hss}D}.~S. 
{\leavevmode\lower.6ex\hbox to 0pt{\hskip-.23ex \accent"16\hss}D}uki{\'c}.
Noether's theorem for optimum control systems.
{\em Internat. J. Control (1)}, 18:667--672, 1973. 
\textsf{Zbl} {0281.49009} \textsf{MR} {49:5979}

\bibitem{MR83k:49054}
A.~van~der Schaft.
Symmetries and conservation laws for {H}amiltonian systems with
inputs and outputs: a generalization of {N}oether's theorem.
{\em Systems Control Lett.}, 1(2):108--115, 1981/82. 
\textsf{Zbl} {0482.93038} \textsf{MR} {83k:49054}

\bibitem{MR2001c:49001}
R.~Vinter.
{\em Optimal control}. Birkh\"auser Boston Inc., Boston, MA, 2000. 
\textsf{Zbl} {0952.49001} \textsf{MR} {2001c:49001}

\bibitem{MR41:4337}
L.~C. Young.
{\em Lectures on the calculus of variations and optimal control theory}.
W. B. Saunders Co., Philadelphia, 1969. 
\textsf{Zbl} {0177.37801} \textsf{MR} {41:4337}

\end{thebibliography}
\end{document}